\documentclass[journal]{IEEEtran}
\usepackage{amsbsy}
\usepackage{floatflt} 

\usepackage{tikz}
\usetikzlibrary{arrows.meta}
\usepackage{amsmath}
\usepackage{amssymb}
\usepackage{times}
\usepackage{graphics}
\usepackage{graphicx}
\usepackage{xspace}
\usepackage{paralist} 
\usepackage{setspace} 
\usepackage{xypic}
\xyoption{curve}
\usepackage{latexsym}
\usepackage{theorem}
\usepackage{ifthen}
\usepackage{caption}
\usepackage{subfigure}
\usepackage{booktabs}
\usepackage{algorithm}
\usepackage{algorithmic}
\usepackage{color}

\usepackage{tabularx,ragged2e}
\newcolumntype{C}{>{\Centering\arraybackslash}X}
\usepackage{nidanfloat}
\newcommand{\fedavg}{\texttt{FedAvg}\xspace}
\newcommand{\fedprox}{\texttt{FedProx}\xspace}

\ifCLASSINFOpdf
\else
\fi

{\theoremheaderfont{\it} \theorembodyfont{\rmfamily}

}

{\theoremheaderfont{\it} \theorembodyfont{\rmfamily}

}

{\theoremheaderfont{\it} \theorembodyfont{\rmfamily}

}

\begin{document}
\title{Primal-dual methods for large-scale and distributed convex optimization and data analytics}

\author{Du$\check{\mbox{s}}$an Jakoveti\'c, Dragana Bajovi\'c, Jo\~ao Xavier, Jos\'e M. F. Moura
\thanks{D. Jakoveti\'c is with Faculty of Sciences, University of Novi Sad, Novi Sad, Serbia. D. Bajovi\'c is with Faculty of Technical Sciences, University of Novi Sad, Novi Sad, Serbia.
 J. Xavier is with Instituto Superior Te\'cnico, Universidade de Lisboa, Lisbon, Portugal, and also with the Institute for Systems and Robotics, Laboratory for Robotics and Engineering Systems.
 J. M. F. Moura is with the Department of Electrical and Computer Engineering,
Carnegie Mellon University,
Pittsburgh, PA, USA. The research of D. Jakovetic is supported by Ministry of Education, Science and Technological Development,
Republic of Serbia, grant no. 174030.
The research of D. Bajovic is supported by the Serbian Ministry
of Education, Science, and Technological Development, Grant no. TR32035.
 The work of J. Xavier is supported
in part by the Fundacao para a Ciencia e Tecnologia, Portugal, under Project
UID/EEA/50009/2019 and Project HARMONY PTDC/EEI-AUT/31411/2017  (funded by Portugal 2020 through FCT,
Portugal, under Contract AAC n 2/SAICT/2017 -- 031411. IST-ID is funded by POR Lisboa -- LISBOA-01-0145-FEDER-031411).
The work of J. M. F. Moura is partially supported by NSF grants CCF \#1513936 and
CPS \#1837607.
Author's e-mails: djakovet@uns.ac.rs; dbajovic@uns.ac.rs;  jxavier@isr.ist.utl.pt; moura@ece.cmu.edu.}}

%
%


\maketitle

\begin{abstract}
The augmented Lagrangian method (ALM) is a classical optimization tool
that solves a given ``difficult'' (constrained) problem via finding solutions of a sequence
of ``easier'' (often unconstrained) sub-problems with respect to the original (primal) variable,
 wherein constraints satisfaction is controlled via the so-called dual variables.
   ALM is highly flexible with respect to how primal sub-problems
   can be solved, giving rise to a plethora of different primal-dual methods.
    The powerful ALM mechanism has recently proved to be
    very successful in various large scale and distributed applications.
    In addition, several significant advances have appeared,
     primarily on precise complexity results with respect to
     computational and communication costs in the presence of inexact updates
     and design and analysis of novel optimal methods for distributed consensus optimization.
      We provide a tutorial-style introduction to ALM and its variants
      for solving convex optimization problems in large scale and distributed settings.
      We describe control-theoretic tools for the algorithms' analysis and design,
      survey recent results, and provide novel insights in the context of two emerging applications:
      federated learning and distributed energy trading.
\end{abstract}


\maketitle \thispagestyle{empty} \maketitle
%
%
%
%

\vspace{-4mm}

\section{Introduction}
\label{section-introduction}

Primal-dual optimization methods
have a long history, e.g.,
\cite{BertsekasBook,Hestenes1969,PowellClassical,ArrowHurwitzUzawa}.
  Seminal works include, e.g., the Arrow-Hurwitz-Uzawa
 primal-dual dynamics~\cite{ArrowHurwitzUzawa}, and the proposal of
 augmented Lagrangians and multiplier methods, e.g.,~\cite{Hestenes1969,PowellClassical}; see also~\cite{HaarhoffBuys1970}.

The augmented Lagrangian (multiplier) method (ALM) is a general-purpose iterative
solver  for \emph{constrained} optimization problems. The underlying mechanism
 translates the task
of solving the original ``difficult" optimization problem with respect to
the (primal) variable~$x$ into solving a sequence
of ``easier'' problems (often unconstrained ones) with respect to~$x$, wherein the
influence of the original ``hard'' constraints is controlled
via the so-called dual variables. The ALM mechanism is beneficial in many scenarios;
for example, it allows to naturally decompose an original large-scale problem into a
set of smaller sub-problems that can then be solved in parallel.
 ALM can handle very generic problems, including non-smoothness of the objective function and generic problem constraints
  and has strong convergence guarantees, e.g.,~\cite{BertsekasBook,ALclassicalBertsekas,ALclassicalAPP3Birgin,ALclassicalAPP2Birgin,ALclassicalAPPpowerSystems}.

The past decade shows a strong renewed interest in primal-dual and
(augmented) Lagrangian methods. This is because they are, by design, amenable to large-scale and distributed
optimization. Indeed, very good performance has been achieved on various Big Data analytics models, including, e.g.,
sparse regression~\cite{ADMMbigDataV1,ADMMbigDataV2,Joao-Mota-2,Joao-Mota-3,cooperative-convex}, large-scale model predictive control~\cite{NecoaraMPC2}, and low rank tensor recovery~\cite{ALnewAPPnetworkAnomalies}, and also on various modern real-world applications, such as  state estimation for the smart grid, e.g.,~\cite{ALnewAPPstateEstimationSmartGrid}, community detection in social networks, e.g.,~\cite{ALnewAPPcommunityDetection}, image filtering, inpainting and demosaicing, e.g.,~\cite{ALnewAPP2wotaoYinTotalVariation,ALnewAPP3TotalVariation,ALnewAPPMarioFigueiredo},
factorization problems in computer vision, e.g., \cite{ALnewAPPbalmXavier}, etc.

In this paper, we provide an overview of
primal-dual methods for large-scale and distributed convex optimization.
  We provide a gentle introduction to the topic,
 followed by an overview of recent results, also suggesting novel insights and
 potential novel applications. In more detail,
 we focus on the following key aspects.

\begin{itemize}
\item \textbf{Approximate updates}.
Several recent results consider the scenarios when primal and/or dual
variables are updated in an inexact way. The works
then study how the inexactness reflects convergence and convergence rate properties of the methods.

\item \textbf{Distributed optimization}.
Primal-dual methods have proved to be very successful in distributed
consensus optimization, e.g., \cite{nedic_T-AC}. Therein, a set of nodes connected in a generic network
collaborate to solve an optimization problem where the overall objective function
is a sum of the components known only locally by individual nodes.
 In this setting, after an appropriate reformulation of the problem of interest, application
 of primal-dual methods leads to efficient and often to optimal methods. Furthermore,
 several existing methods
 that have been derived from different perspectives
 have been recently shown to admit primal-dual interpretations,
 which in turn opens up possibilities for further method improvements.

\item \textbf{Complexity results}. While traditional studies of primal-dual and ALM
focus on convergence rate guarantees with respect to the (outer) iteration count,
more recent works focus on establishing complexity results with respect to
more fine-grained measures of communication and computational costs,
such as number of gradient evaluations, number of inner iterations, etc.
 These results give additional insights into the performance of the methods
 with respect to the actual communication and/or computational costs.

\item \textbf{Control-theoretic analysis and interpretations}. It has been recently shown
that primal-dual and ALM admit interpretations from
the control theory perspective. This allows to
interpret, analyze and design the methods with classical control-theoretic tools.
This makes the methods accessible to a broader audience than before but also
allows for improving the methods.

\item \textbf{Applications}.
We demonstrate how an instance of ALM-type methods, namely
the parallel direction method of multipliers~\cite{PDMMTomLuo}, can be applied to
the emerging concept of federated learning.
 We capitalize on recent results on inexact ALM
to show how some recently proposed energy trading methods
 can be modified to provably work,
with (inner) iteration complexity certifications,
under more general energy generation and energy transmission
cost models.
\end{itemize}

We now contrast our paper with other tutorial and survey papers
 on primal-dual or decentralized methods, e.g., \cite{BoydADMM,NedicRabbatProceedingsIEEE}.
 Reference~\cite{BoydADMM} is concerned with ADMM, a specific
 primal-dual method, and it considers master-worker (star network)
  architectures, i.e., it is not concerned with fully distributed
  network architectures. In contrast,
  we overview here both the methods amenable to master-worker
   and to fully distributed topologies.
  Reference~\cite{NedicRabbatProceedingsIEEE}
   provides an overview of
   distributed optimization
   methods, focusing on complexity, communication-computational tradeoffs,
   and network scaling. In doing so, the paper mainly
   treats \emph{primal} distributed (sub)gradient methods.
   In contrast, we focus here on
   \emph{primal-dual} methods for
   large scale and distributed convex optimization.

   In summary, our motivation
   for providing this overview is to present to a broader readership
   the following.
    1) Primal-dual methods exhibit a high degree of
    flexibility, in terms of the underlying network
    topology to which they can be adapted (e.g., master-worker, fully
    distributed); 2) they exhibit a high degree of flexibility
    with respect to
    the ways (inexact) primal variable updates can be performed;
    3) they are amenable to control-theoretic analysis and design;
     4) several state of the art existing methods that were originally proposed
     from a different perspective (see, e.g., \cite{Shi2015,harnessing,SmallGainNedicUncoord}) can be shown to be instances
     of primal-dual methods; and 5) primal-dual methods, aproprietly designed,
     exhibit good, and often optimal, computation-communication complexity.

\textbf{Paper outline}. Section~\ref{section-master-worker}
provides preliminaries on ALM and related methods and presents how they can be applied
in distributed optimization.
In particular, Subsection~{II-A} gives preliminaries on ALM. Then, we consider
more traditional master-worker architectures (Subsection~{II-B}),
followed by a treatment of fully distributed architectures (Subsection~{II-C}).
 Section~{III} provides control-theoretic analysis and insights into
 distributed ALM methods. Section~{IV} provides a review of recent
 results on the topic. Section~V considers two applications --
  one for master-worker architectures (federated learning, Subsection~V-A),
  and one for fully distributed architectures (energy trading in microgrids, Subsection~V-B).
  Finally, we conclude in Section~{VI}.

\vspace{-4mm}

\section{Primal-dual methods for distributed optimization}
\label{section-master-worker}
Subsection~{II-A} briefly introduces ALM and some related methods;
Subsection~{II-B} shows how they can be applied on master-worker computing architectures; and
Subsection~{II-C} describes how they can be applied on fully distributed architectures.

\subsection{Preliminaries}
\label{subsection-preliminaries-AL}

We now introduce the classical ALM. To illustrate, consider an optimization problem with linear constraints,
\begin{equation} \begin{array}[t]{ll} \underset{x \in {{\mathbb R}}^n}{\text{minimize}} & f(x) \\ \text{subject to} & A x = 0, \end{array} \label{eq:ALexample} \end{equation}
where $f\colon {{\mathbb R}}^n \rightarrow {{\mathbb R}}$ is the given objective function and $A$ is a given matrix in ${{\mathbb R}}^{m\times n}$.
Starting from an arbitrary point $\left( {x^{(0)}}, {\lambda^{(0)}} \right) \in {{\mathbb R}}^n\times {{\mathbb R}}^m$, the ALM generates a sequence of points $\left\{ \left( {x^{(k)}}, {\lambda^{(k)}} \right) \right\}_{k \geq 1}$
as
\begin{eqnarray}
{x^{(k+1)}} & = & \arg\min_{x \in {{\mathbb R}}^n}\, {\mathcal L}\left( x , {\lambda^{(k)}} \right) \label{eq:ALprimalupdate} \\ \lambda^{(k+1)} & = & {\lambda^{(k)}} + \rho A {x^{(k+1)}}. \label{eq:ALdualupdate}
\end{eqnarray}
Here, ${x^{(k)}}$ and ${\lambda^{(k)}}$ are primal and dual variables, respectively.
 Further, the function ${\mathcal L}( x, \lambda )$, which is defined as \begin{equation} {\mathcal L}( x, \lambda) = f(x) + \lambda{^\top} A x + \rho \left\| A x \right\|^2 / 2 \label{eq:ALfunction} \end{equation}
and is called the augmented Lagrangian of problem~\eqref{eq:ALexample}, determines the  functions being successively minimized without constraints in the ALM update~\eqref{eq:ALprimalupdate}. The constant $\rho$ is a previously chosen positive number.

The convergence properties of the ALM have been studied in much detail, both for the convex and nonconvex settings;
reference~\cite{Bertsekas1982} is standard, while~\cite{Nedic2008} lays out a more recent geometric framework, and~\cite{FernandezSolodov2012} proves convergence under unusually weak constraint qualifications.

The convergence properties of the ALM are particularly appealing in the convex setting. Indeed,
 the landmark paper~\cite{Rockafellar1976} links the ALM to the proximal point method for monotone operators,\footnote{
The proximal point method for unconstrained minimization of
a convex function~$\phi: \,{\mathbb R}^n \rightarrow \mathbb R$
 works as follows: $x^{(k+1)} = \mathrm{arg\,min}_{x \in {{\mathbb R}^n}}
 \{\phi(x)+\frac{1}{2 \,\alpha}\|x-x^{(k)}\|^2\}$,
 where $x^{(k)}$ is the solution estimate at iteration~$k$,
 and $\alpha$ is a positive parameter. The method converges
 at rate~$1/k$ for convex functions; see, e.g.,~\cite{Guler}.}
 showing that the ALM is just the proximal point method applied to the dual problem of~\eqref{eq:ALexample},
\begin{equation} \begin{array}[t]{ll} \underset{\lambda \in {{\mathbb R}}^m}{\text{maximize}} & \mathcal{D}(\lambda), \end{array} \label{eq:dualALexample} \end{equation}
where $\mathcal{D}$ is the dual function $\mathcal{D}( \lambda ) = \inf\{ f(x) + \lambda{^\top} A x \colon x \in {{\mathbb R}}^n \}$. This link transfers at once the strong convergence properties of the proximal point method to the ALM, a typical convergence result being that the sequence $\left\{ {\lambda^{(k)}} \right\}_{k \geq 1}$ converges to a solution $\lambda^\star$  (assuming it exists) of the dual problem~\eqref{eq:dualALexample}, while the sequence $\left\{ {x^{(k)}} \right\}_{k \geq 1}$ converges to a solution $x^*$ of the primal problem~\eqref{eq:ALexample} (assuming, for example, that $f$ is strongly convex).

An important feature of ALM is that the update in \eqref{eq:ALprimalupdate} can be done in an inexact fashion,
giving the method a high degree of flexibility, but also giving rise to various
method variants. Specifically, if instead of \eqref{eq:ALprimalupdate} we carry out a single gradient step
 on function $\mathcal{L}$ with respect
 to variable $x$, we recover a method closely related to the
 Arrow-Hurwitz-Uzawa (AHU) dynamics (saddle point method)~\cite{ArrowHurwitzUzawa}:
\begin{eqnarray}
{x^{(k+1)}} & = &
x^{(k)} - \alpha\,\, \nabla_x{\mathcal L}\left( x^{(k)} ,
{\lambda^{(k)}} \right) \label{eq:ALprimalupdate-AHU-Sec-2} \\
\lambda^{(k+1)} & = & {\lambda^{(k)}} + \rho A {x^{(k+1)}}, \label{eq:ALdualupdate-AHU-Sec-2}
\end{eqnarray}
where constant $\alpha>0$ is a primal step-size.\footnote{
Note that \eqref{eq:ALdualupdate-AHU-Sec-2}
can be written as
$\lambda^{(k+1)}  =  {\lambda^{(k)}} + \rho \,\nabla_{\lambda}{\mathcal L}\left( x^{(k+1)} ,
{\lambda^{(k)}} \right)$, i.e.,
\eqref{eq:ALdualupdate-AHU-Sec-2} is a dual gradient ascent step.
Clearly, the same applies to~(3). In addition, in turns out
that~(3) corresponds to an ascent step according to
the gradient of the \emph{dual function} $\mathcal{D}$ in~(5).
}
 Similarly, when the optimization variable is
 partitioned into $2$ blocks $x=((x_1)^\top,(x_2)^\top)^\top$,
 and when the objective function and the constraint can be written in the following form:
$
 f(x) =  f_1(x_1) + f_2(x_2)$,  $A\,x  =  A_1\,x_1+A_2\,x_2=0,$
  then updating $x_1$ and $x_2$ sequentially in a Gauss-Seidel fashion, followed by
 the dual update \eqref{eq:ALdualupdate}, gives rise
 to the celebrated alternating direction method of multipliers (ADMM), e.g., \cite{BoydADMM}.
   See also, e.g., \cite{NecoaraAUFLagr,GeorgeLanAL}, for analysis
   of inexact augmented Lagrangian methods.

\subsection{Dual decomposition framework: Master-worker architectures}
\label{subsection-conventional-decomposition}
The conventional (dual) decomposition framework assumes a
computational infrastructure with
$N$ computing nodes (workers), each of which can communicate
with a master node.
 The relevant optimization problem
 then takes the form \eqref{eq:ALexample},
 with the objective function and the constraint given by:
 $
 f(x) =  f_1(x_1) + f_2(x_2)+...+f_N(x_N)$;
  $A\,x  =  A_1\,x_1+A_2\,x_2+...+A_N\,x_N=0.
  $
 Here, we assume that the optimization variable $x = (x_1^\top,...,x_N^\top)^\top$
is partitioned into $N$ blocks $x_i$'s, $i=1,...,N$, where
the $i$-th block is assigned to the $i$-th worker.
Similarly, the objective function summand $f_i(x_i)$ is assigned to worker $i$.
Matrix $A$  is partitioned into blocks $A_i$'s, where all the $A_i$'s are available to
the master and to all workers.
  We assume that
 worker~$i$ is responsible for updating the primal variable block $x_i$,
  while the master is responsible for updating the dual variable.

The augmented Lagrangian function $\mathcal{L}$
can now be written as follows:
\begin{equation}
\label{eqn-AL-dual-decomposition}
\mathcal{L}(x,\lambda) = \sum_{i=1}^N f_i(x_i) + \lambda^\top \left( \sum_{i=1}^N A_i\,x_i\right)
+\frac{\rho}{2} \left\| \sum_{i=1}^N A_i\,x_i \right\|^2,
\end{equation}
and the ALM proceeds as in \eqref{eq:ALprimalupdate}--\eqref{eq:ALdualupdate}.
%
%
%
 Due to the quadratic term in~\eqref{eqn-AL-dual-decomposition},
the update~\eqref{eq:ALprimalupdate} cannot be executed in parallel,
in the sense that each of the workers $i$ updates $x_i$ in parallel
with other workers.
 To overcome this issue,
several algorithms that parallelize~\eqref{eq:ALprimalupdate}
have been proposed, including the diagonal quadratic approximation, e.g.,
\cite{Ruszcynski}, and
more recently~\cite{ZavlanosMethod}.
As presented in~\cite{PDMMTomLuo}, an appealing attempt is to update the primal updates in a Jacobi fashion:
\begin{align}
 x_i^{(k+1)} &= \mathrm{arg\,min}_{x_i}~ \mathcal{L}\left(x_i, x_{j\neq i}^{(k)}, \lambda^{(k)} \right)~, \label{eqn-PJ-DD-1}\\
 \lambda^{(k+1)} &=  \lambda^{(k)} + \rho A\,x^{(k+1)}. \label{eqn-PJ-DD-2}
\end{align}
That is, each worker $i$ optimizes function $\mathcal{L}$
with respect to $x_i$ in parallel,
while the values of the other workers'
variables are fixed to their
value before the current parallel update, i.e.,
$x_j$ is fixed to $x_j^{(k)}$, $j \neq i$.
However,
 it is known, e.g., \cite{PDMMTomLuo},
  that the method in \eqref{eqn-PJ-DD-1}--\eqref{eqn-PJ-DD-2}
  may not converge. An interesting recent method,
  dubbed parallel direction method of multipliers (PDMM)~\cite{PDMMTomLuo},
  overcomes this problem and works as follows.
 At each iteration $k$, PDMM
 first randomly selects $K$ out of $N$ variable blocks,
 where the set of selected blocks is denoted by $\mathbb{I}_k$.
 Then, the following update is carried out:
 {\allowdisplaybreaks{
\begin{align}
x_{i}^{(k+1)} &= \mathrm{arg\,min}_{x_{i}}~ \left\{ \mathcal{L}
(x_i, x_{j\neq i}^{(k)}, \hat{\lambda}^{(k)}) \right. \nonumber \\
&+ \left.
\eta_{i}^{(k)} B_{i}(x_{i}, x_{i}^{(k)}) \right\},\,\,\,
i \in\mathbb{I}_k, \label{eqn-rpdmm_xi} \\
 \lambda^{(k+1)} &=  \lambda^{(k)} + \Theta \,\rho A\,x^{(k+1)}~,\label{eqn-rpdmm_y}  \\
  \hat{\lambda}^{(k+1)} & = \lambda^{(k+1)} - \Theta^\prime \rho A\,x^{(k+1)} , \label{eqn-rpdmm_yhat}
\end{align}}}
Here, $\Theta = \mathrm{Diag}\left( \tau_1,...,\tau_N\right)$ is a diagonal
matrix of positive weights $\tau_i$; $\Theta^\prime = \mathrm{Diag}\left( \nu_1,...,\nu_N\right)$,
$\nu_i \in [0,1)$;  $\eta_{i}^{(k)} > 0$; and   $B_i(\cdot,\cdot)$ is a Bregman divergence.\footnote{
Given a continuously differentiable
convex function $\Psi: {\mathbb R}^n \rightarrow \mathbb R$,
the Bregman divergence $B: {\mathbb R}^n \times {\mathbb R}^n \rightarrow \mathbb R$ is defined as:
 $B(u,v) = \Psi(u) - \Psi(v) - \left(\nabla \Psi(v)\right)^\top \,(u-v)$.
  For example, for $\Psi(u) = \|u\|^2$, we get as Bregman divergence
  the Euclidean distance $B(u,v)=\|u-v\|^2$;
   for $\Psi(u) = \frac{1}{2}u^\top A u$,
   we get $B(u,v) = \frac{1}{2}(u-v)^\top A (u-v),$ where $A$ is a positive definite matrix.
 }

PDMM differs from \eqref{eqn-PJ-DD-1}--\eqref{eqn-PJ-DD-2} in three main aspects.
First, PDMM updates only a subset of primal variables $x_i$'s at each iteration.
Second, the primal variable update involves a Bregman divergence, which
makes the primal variable trajectory ``smoother.'' Thirdly, PDMM
 introduces a ``backward'' step~\eqref{eqn-rpdmm_yhat} in the dual
 variable update to further smooth the primal-dual variable trajectory.
 We refer the reader to \cite{PDMMTomLuo}, equation~(9),
 regarding the motivation for the backward step.
 Interestingly, when the constants $\eta_i^{(k)}$'s
 are sufficiently large, the backward step~\eqref{eqn-rpdmm_yhat}
 is not needed, and variables $\lambda^{(k)}$ and $\hat{\lambda}^{(k)}$ coincide.
 We will further exploit PDMM in the context of federated learning in Section~V-A.

\subsection{Primal-dual methods for consensus optimization}
\label{subsection-fully-distributed}

We next consider distributed consensus optimization, e.g., \cite{nedic_T-AC,SoummyaGossipIEEEproc,AngeliaIEEEProc},
wherein a set of $N$ agents collaborate to solve the following problem:\footnote{For notational simplicity,
we let the optimization variable in \eqref{eq:dopt} be a scalar one;
it can be easily extended to vector variables in ${\mathbb R}^d$, but the
notation becomes more cumbersome, involving extensive use
of Kronecker products; see, e.g., \cite{Shi2015}.
}
\begin{equation} \begin{array}[t]{ll} \underset{x \in {{\mathbb R}}}{\text{minimize}} & f_1(x)  + \cdots + f_N(x), \end{array} \label{eq:dopt} \end{equation}
Here, the convex function  $f_i\colon {{\mathbb R}} \rightarrow {{\mathbb R}}$ is known only at agent~$i$. To achieve this goal, each agent can communicate with a handful of other agents through an underlying sparse communication network, assumed undirected and connected. Note that this key object---the communication network---is absent from the formulation~\eqref{eq:dopt}. A common ploy  to bring it to view is to reformulate~\eqref{eq:dopt} as
\begin{equation} \begin{array}[t]{ll} \underset{x = ( x_1, \ldots, x_N) \in {{\mathbb R}}^N}{\text{minimize}} & f_1(x_1)
 + \cdots + f_N(x_N) \\ \text{subject to} & \rho\,L^{1/2} x = 0, \end{array} \label{eq:dopt2} \end{equation}
where
$\rho$ is a positive parameter and $L$ is the weighted graph Laplacian
of the communication network (e.g., see~\cite{Mokhtari2016,UnificationAndGeneral}), assumed to be connected.
The constraint, which involves the matrix square-root of the semidefinite positive matrix $L$,
ensures that all components of $x$ are equal, that is, $x_1 = \cdots = x_N$,
thereby securing the equivalence between~\eqref{eq:dopt} and~\eqref{eq:dopt2}.
 Let $f\colon {{\mathbb R}}^N \rightarrow {{\mathbb R}}$ be defined as
$f(x_1, \ldots, x_N) = f_1( x_1) + \cdots + f_N( x_N ).$
The augmented Lagrangian and the ALM~\eqref{eq:ALprimalupdate}--\eqref{eq:ALdualupdate} to solve problem~\eqref{eq:dopt}
 are then, respectively, given by:
 \begin{eqnarray}
 \mathcal{L}(x,\lambda) &=& f(x) + \rho\,\lambda^\top L^{1/2}x
 + \frac{\rho}{2}x^\top L x \label{eqn-AL-consensus-new}\\
  x^{(k+1)} &=& \mathrm{arg\,min}_{x \in {{\mathbb R}}^N}~\mathcal{L}\left(x, \lambda^{(k)} \right)
   \label{eqn-AL-generic-ALM-primal} \\
 \lambda^{(k+1)} &=&  \lambda^{(k)} + {\rho} \,\nabla_{\lambda} \mathcal{L}(x^{(k+1)},\lambda^{(k)}),
 \label{eqn-AL-generic-ALM-dual}
 \end{eqnarray}
where $\nabla_{\lambda}$ denotes the partial derivative with respect to variable~$\lambda$.
The algorithm above further simplifies to the following:
 \begin{eqnarray}
{x^{(k+1)}} & = & \arg\min_{x \in {{\mathbb R}}^N} f(x) + x^\top  \eta^{(k)} +
 \frac{\rho}{2} x{^\top} L x  \label{eq:ALp} \\ {\eta^{(k+1)}} & = & {\eta^{(k)}} + \rho\, L {x^{(k+1)}}, \label{eq:ALd}
\end{eqnarray}
where the original variable ${\lambda^{(k)}}$ is swapped for the novel one ${\eta^{(k)}} =
{\rho} L^{1/2} {\lambda^{(k)}}$
 to pull off the problematic matrix $L^{1/2}$ from the formulation
 (as its sparsity does not mirror the sparsity of the communication network).\footnote{Note that, in this way, the sparsity-incompatible matrix $L^{1/2}$ is only used to construct
the problem reformulation in~\eqref{eq:dopt2}, while it does
 not play a role in the actual computations in~\eqref{eq:ALd}.}

In this formulation, agent~$i$ manages the $i$-th components of ${x^{(k)}}$ and ${\eta^{(k)}}$.
Also, the update \eqref{eq:ALd} is implementable in a fully distributed way,
as matrix $L$ respects the underlying graph sparsity. However, similarly to
the dual decomposition framework scenario in Subsection~{II}-A,
the quadratic term $x{^\top} L x$ in the primal update~\eqref{eq:ALp}
couples all the agents, spoils the otherwise nice separable structure of the primal update,
and blocks a directed implementation of this method in a distributed setting.

This grim state of affairs, however, is easy to lift because the canonical ALM can withstand several
changes without losing its main convergence properties.
That is, the canonical ALM~\eqref{eq:ALp}--\eqref{eq:ALd}
can be used as a starting point to inspire more distributed-friendly algorithms.

 In fact, there is a very high degree of flexibility
 of how one can perform the primal update~\eqref{eq:ALp}
 and the dual update~\eqref{eq:ALd}; this gives rise to a plethora of different
 distributed methods for consensus optimization.
  We mention here a few examples while
  a more complete overview of the possibilities is provided
  in Subsection~\ref{subsection-overview-consensus} and Table~I.
   For example, the ADAL distributed method~\cite{Chatzipanagiotis2015} changes problem~\eqref{eq:ALp} to one that
separates across agents and adjusts the dual update~\eqref{eq:ALd} accordingly.
In a different direction, reference~\cite{LinearRateAL} changes the exact minimization of~\eqref{eq:ALp}
to an inexact one, showing that inexact solves of~\eqref{eq:ALp}, when properly controlled,
not only retain convergence but also are amenable to distributed computation (say, by tackling~\eqref{eq:ALp}
with some iterations of the gradient method, each of which is directly distributed).

An even more aggressive change to~\eqref{eq:ALp} is to carry out but one gradient step from the current iterate,
i.e., carrying out an AHU-type algorithm. This also gives an advantage of
avoiding ``double-looped'' methods, where primal and dual updates are carried out
at different time scales. The AHU-type
algorithm blueprint leaves us near the structure of efficient distributed first-order primal-dual methods
 such as EXTRA~\cite{Shi2015}, DIGing~\cite{SmallGainNedicTimeVar}, and their subsequent generalization in~\cite{UnificationAndGeneral}.
 More precisely, while EXTRA and DIGing have been developed from a perspective
different than primal-dual, they are in fact instances of primal dual methods
 as first shown in~\cite{SmallGainNedicTimeVar}.

To illustrate, we consider EXTRA proposed in~\cite{Shi2015},
  a method developed from the standpoint
of correcting the fixed point
equation of the standard distributed gradient method~\cite{nedic_T-AC}:
\begin{eqnarray*}
\label{eqn-extra}
x^{(1)} &=& \hspace{-2mm} {W}\,x^{(0)} - \alpha \,\nabla f(x^{(0)})\\
\label{eqn-extra-222}
x^{(k+1)} &=&  \hspace{-2mm}
2\, {W}\, x^{(k)} -
\alpha\,\nabla f(x^{(k)})- {W}\,x^{{(k-1)}}  \\
&\,&+\,
\alpha\,\nabla f(x^{(k-1)}),\,\,\,k=1,2,... ,\nonumber
\end{eqnarray*}
where $W=I-L$ is the weight matrix and $\alpha>0$ is a step-size. The matrix $W$ is doubly stochastic and  respects the sparsity of the underlying graph.
  For the AL function in~\eqref{eqn-AL-consensus-new}, it can be shown
that EXTRA is equivalent to the following method, with $\alpha=1/\rho$~\cite{Shi2015}; see also~\cite{UnificationAndGeneral}:
 \begin{eqnarray}
 \label{eqn-saddle-point-Extra}
 x^{(k+1)} &=& x^{(k)} - \alpha\, \nabla_x \mathcal{L}(x^{(k)},\lambda^{(k)})\\
 \lambda^{(k+1)} &=&
 \lambda^{(k)} +  \alpha\,\nabla_{\lambda} \mathcal{L}(x^{(k+1)},\lambda^{(k)}).
 \end{eqnarray}

On the other hand, the methods in~\cite{harnessing,SmallGainNedicTimeVar}, see also~\cite{SmallGainNedicUncoord,SmallGainKin1},
 utilize a gradient tracking principle and work as follows:
 \begin{eqnarray*}
\label{eqn-harnessing-dis-opt-sect}
x^{(k+1)} &=& {W}\,x^{(k)} - \alpha\,s^{(k)}\\
\label{eqn-harnessing-primeprime}
s^{(k+1)} &=& {W}\,s^{(k)} + \nabla f(x^{(k+1)}) - \nabla f(x^{(k)}),
\end{eqnarray*}
where $s^{(k)}$
 is an auxiliary variable
 that tracks the global average of the $f_i$'s
 individual gradients along iterations.
 The  method also
 accounts for a primal dual interpretation, as shown in~\cite{SmallGainNedicTimeVar,SmallGainNedicUncoord};
 see also~\cite{UnificationAndGeneral}. For example, following~\cite{UnificationAndGeneral},
  update rule~\eqref{eqn-harnessing-dis-opt-sect}--\eqref{eqn-harnessing-primeprime}
  can be rewritten as follows (here, $\alpha=1/\rho$):
 \begin{eqnarray}
 \label{eqn-saddle-point-Extra}
 x^{(k+1)} &=& x^{(k)} - \alpha\, \nabla_x \mathcal{L}(x^{(k)},\lambda^{(k)})\\
 \lambda^{(k+1)} &=&
 \lambda^{(k)} +  \alpha\,\nabla_{\lambda} \mathcal{L}(x^{(k+1)},\mu^{(k)})\\
 &-& \alpha\,W\,\nabla_{\lambda} \mathcal{L}(x^{(k)},\mu^{(k)}).
 \end{eqnarray}

The primal-dual methodology allows not only
to re-interpret many existing first order methods, but it also
allows for derivation of some novel methods.  For example,
reference~\cite{UnificationAndGeneral} proposes a method that
introduces an additional tuning parameter matrix $B$ that has the sparsity
pattern reflecting the network.   Reference~\cite{UnificationAndGeneral} shows that
specific choices of matrix $B$ recover the EXTRA and DigING
methods. The author derives a primal-dual error recursion for
the proposed method, which allows to tune parameter $B$ for
an improved algorithm performance.


\vspace{-4mm}

\section{Analysis and insights}
\label{section-analysis-insights}

We now highlight some of the ideas spawned by a recent thread, fertile and ongoing,  of works that look at  distributed optimization from a control-theoretic viewpoint.
The control toolset enriches the analysis and design of optimization algorithms with many exciting insights and proof techniques, too many to cover in detail here; with the beginner in mind, we illustrate the use of the basic tool of the LaSalle invariance principle as a proof technique and the interpretation of the ALM as a basic proportional-integral controller.

\subsection{The LaSalle invariance principle}
\label{subsection-laSalle}
 The LaSalle invariance principle is a standard tool used in control to analyze the behavior of nonlinear dynamical systems. The LaSalle invariance principle is covered, e.g.,  in the book~\cite{Khalil2002}, which focuses on  continuous-time systems; discrete-time systems are studied in greater depth in the recent tutorial~\cite{Mei2019}.
Roughly speaking, LaSalle's principle characterizes the fate of trajectories of a given nonlinear dynamical system, say \begin{equation} \dot x(t) = \Phi(x (t)), \label{eq:dynsys} \end{equation} where $\Phi\colon {{\mathbb R}}^N \rightarrow {{\mathbb R}}^N$ is assumed here, for simplicity, to guarantee trajectories for any initial condition $x(0)$ and for all $t \in {{\mathbb R}}$ (that is, no trajectory escapes in finite-time). Suppose $V\colon {{\mathbb R}}^N \rightarrow {{\mathbb R}}$ is a Lyapunov-like function for the given dynamical system, in the sense that $V$ is smooth, has bounded sublevel sets ($\{ V \leq c \} \subset {{\mathbb R}}^N$ is compact for any $c$), and does not increase along trajectories of the dynamical system ($\dot V(x(t)) \leq 0$ for all $t \in {{\mathbb R}}$, whenever $x\colon {{\mathbb R}}^N \rightarrow {{\mathbb R}}$ is a solution of~\eqref{eq:dynsys}). Then, the LaSalle invariance principle asserts that each trajectory converges to the largest invariant subset $M$ of the set $E = \{ x \in {{\mathbb R}}^N \colon \nabla V \cdot \Phi(x) = 0 \}$ (usually, denoted $E = \{ \dot V = 0 \}$), a set $A$ being called invariant if trajectories starting in $A$ are fully contained in $A$.

As an illustration, we show how LaSalle's principle offers a quick proof of convergence for the distributed algorithm in~\cite{Wang2010}:
\begin{eqnarray}
\dot x(t) & = & -\nabla f( x(t) ) - L \eta(t) - L x(t) \label{eq:LSp} \\ \dot \eta(t) & = & L x(t). \label{eq:LSd}
\end{eqnarray}
Algorithm~\eqref{eq:LSp}-\eqref{eq:LSd} can be interpreted as an AHU dynamics~\eqref{eq:ALprimalupdate-AHU-Sec-2}--\eqref{eq:ALdualupdate-AHU-Sec-2},
or as a version of ALM~\eqref{eq:ALp}-\eqref{eq:ALd} in which the primal update is changed to just a gradient step. The major distinction, however, is that algorithm~\eqref{eq:LSp}-\eqref{eq:LSd} runs in continuous-time. Algorithms in continuous-time cannot be directly implemented in digital setups. Indeed, it is well known that the discrete time and continuous
time algorithms are not equivalent, unless time
scaling or stepsize choices are done properly or implicit discretization is
used; see, e.g., \cite{discretization1, discretization2}.
They are nevertheless worth investigating, if only, as a first step to probe the soundness of a given discrete-time algorithm, because proofs in continuous-time are often much easier and shorter, and can even trigger new insights (e.g., see~\cite{Muehlebach2019} for such a recent example). Continuous-time algorithms for distributed optimization feature, e.g.,  in~\cite{Wang2011,Droge2014,Gharesifard2014,Tang2016,Porco2016}.

Returning to the analysis of algorithm~\eqref{eq:LSp}-\eqref{eq:LSd}, whose purpose is to solve problem~\eqref{eq:dopt2}, we let $x^\star \in {{\mathbb R}}^N$ be a solution of problem~\eqref{eq:dopt2}. We do not assume here such solution is unique but, for simplicity, we assume that the convex function $f$ has a Lipschitiz gradient with constant $C > 0$, that is, the inequality $\left\| \nabla f(x) - \nabla f(y) \right\| \leq C \left\| x - y \right\|$ holds for all $x, y \in {{\mathbb R}}^n$. This property further implies that the gradient is co-coercive, that is, \begin{equation} \left( x - y \right){^\top} \left( \nabla f(x) - \nabla f(y) \right) \geq \frac{\left\| \nabla f(x) - \nabla f(y) \right\|^2}{C}, \label{eq:ccL} \end{equation}
for $x, y \in {{\mathbb R}}^n$.
Although a solution of~\eqref{eq:dopt2} may not be unique, all such solutions share the same gradient:  if $x^\star$ and $y^\star$ solve~\eqref{eq:dopt2}, then $\nabla f\left( x^\star \right) = \nabla f\left( y^\star \right)$. This follows from~\eqref{eq:ccL} and the Karush-Kuhn Tucker conditions for~\eqref{eq:dopt2}, which reveal that the gradient of a  solution $x^\star$ is of the form $-L \eta^\star$, for some $\eta^\star \in {{\mathbb R}}^N$ (and also, of course, $L x^\star = 0$).

Now, consider the function $$V(x, \eta) = \frac{1}{2} \left\| x - x^\star \right\|^2  + \frac{1}{2} \left\| \eta - \eta^\star \right\|^2,$$
which can be readily seen to satisfy
    \begin{equation*}
    \begin{split}
\dot V( x(t), \eta(t) ) = -\left( x(t) - x^\star \right){^\top} \left( \nabla f\left( x(t) \right)  - \nabla f( x^\star ) \right) \\
    - \left( x(t) - x^\star \right){^\top} L \left( x(t) - x^\star \right)
    \end{split}
    \end{equation*}
    for trajectories $( x(t), \eta(t) )$ of the dynamical system~\eqref{eq:LSp}-\eqref{eq:LSd}. From~\eqref{eq:ccL}, we see that $\dot V( x(t), \eta(t) ) \leq 0$ and the LaSalle invariance principle can be applied: it shows that each trajectory converges to the largest invariant set $M$ of the set $E = \{ \dot V = 0 \}$, which, in this example,  amounts to $E = \{ (x, \eta) \colon \nabla f(x) = \nabla f\left( x^\star \right), L x = 0 \}$; so, the point $(x, \eta)$ is in the set $E$ if and only if $x$ solves~\eqref{eq:dopt2}.

The next step is to find the set $M$, which consists of the trajectories fully contained in $E$. Let $\left\{ (x(t), \eta(t) ) \right\}_{t \in {{\mathbb R}}}$ be such a trajectory, that is, $\nabla f( x(t) ) = \nabla f( x^\star )$ and $L x(t) = 0$ for all $t \in {{\mathbb R}}$. These equalities, when plugged in~\eqref{eq:LSp}-\eqref{eq:LSd}, imply that $\eta(t) = \eta$, for some $\eta \in {{\mathbb R}}^n$, and that  $\dot x(t) = - \nabla f( x^\star ) - L  \eta$, for all $t \in {{\mathbb R}}$. Finally, intersecting the last equality with $L x(t) = 0$ for all $t \in {{\mathbb R}}$ gives $\nabla f( x^\star ) = - L \eta$. Therefore, $M = \left\{ (x, \eta) \colon  \nabla f(x) = - L \eta, L x = 0 \right\}$, which is the set of Karush-Kun-Tucker (KKT) points (see, e.g., \cite{BertsekasBook}) for the problem
\begin{equation} \begin{array}[t]{ll} \underset{x = ( x_1, \ldots, x_N) \in {{\mathbb R}}^N}{\text{minimize}} & f_1(x_1)  + \cdots + f_N(x_N) \\ \text{subject to} & L x = 0. \end{array} \label{eq:dopt3} \end{equation}
To conclude, the LaSalle invariance principle says that the distributed algorithm~\eqref{eq:LSp}-\eqref{eq:LSd} navigates any initial point $( x(0), \eta(0) )$ to a set $M$ in which each point is of the form $(x, \eta)$, with $x$ being a solution of~\eqref{eq:dopt2} (or, equivalently,~\eqref{eq:dopt3}) and $\eta$ being a solution  of the dual problem of~\eqref{eq:dopt3}. This conclusion is partially contained in Theorem~3.2 of reference~\cite{Wang2010} (whose proof is omitted therein). Before moving to a new topic in the next section, we point out that other concepts of control, such as dissipativity theory, have also proved useful in establishing convergence of optimization algorithms, e.g., see~\cite{Tang2016,Porco2016}. In addition,
 reference~\cite{Passivity} gives a passivity-based
perspective for a distributed optimization
algorithm studied therein.
\subsection{Proportional-integral controller}
\label{subsection-PI-controller}
\begin{figure}[h]
\begin{center}
\begin{tikzpicture}[scale=0.8]
\draw [thick] (1.5,4) circle [radius=0.25];
\draw[-{Latex[length=2mm]}] (0.5,4)--(1.25,4);
\node [left] at (0.5,4) {$r$};
\draw[-{Latex[length=2mm]}] (1.75,4)--(2.5,4);
\draw[thick] (2.5,3.2) rectangle (4.5,4.8);
\draw[-{Latex[length=2mm]}] (4.5,4)--(5.5,4);
\draw[thick] (5.5,3.2) rectangle (7.5,4.8);
\draw[-{Latex[length=2mm]}] (7.5,4)--(9,4);
\node [right] at (9,4) {$y$};
\node at (3.5,4) {$K$};
\node at (6.5,4) {$G$};
\node [above] at (2,4) {$e$};
\node [above] at (5,4) {$u$};
\draw [-{Latex[length=2mm]},thick] (8,4) -- (8,2) -- (1.5,2) -- (1.5,3.75);
\node [right] at (1.5,3.70) {$-$};
\draw [fill] (8,4) circle [radius=0.1];
\end{tikzpicture}
\end{center}
\caption{A typical control setup: the desired reference signal $r$ is compared to the plant $G$ output signal $y$ to yield the error signal $e$; the controller $K$ monitors the error signal $e$ and reacts by driving the plant with a suitable actuator signal $u$ so as to cancel the error $e$.\label{fig:controlsetup}}
\end{figure}
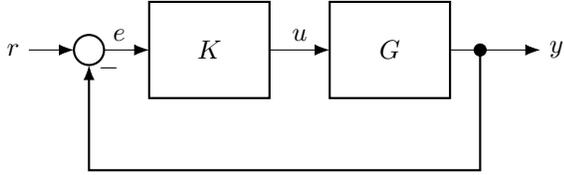

Interestingly, the primal-dual algorithm~\eqref{eq:LSp}-\eqref{eq:LSd} can be interpreted as a widely used controller structure: the proportional-integral controller. Consider the typical control setup depicted in Figure~\ref{fig:controlsetup}, in which the controller system $K$ acts on the plant $G$ so as to make the plant output signal $y$ track a given desired reference signal $r$.

A popular choice for the controller is the proportional-integral (PI) structure, in which the PI controller computes the actuator signal $u(t)$ as a linear combination of a term proportional to the error and a term that integrates the error: $u(t) = \kappa_P e(t) + \kappa_I \int_{-\infty}^t e(s) ds$, where the coefficients $\kappa_P$ and $\kappa_I$ balance the two terms and must be carefully tuned to ensure stability and other performance criteria of the closed-loop system.

  Retracing the ideas presented in~\cite{Wang2011} and~\cite{Droge2014},
  we can cast the primal-dual algorithm \eqref{eq:LSp}-\eqref{eq:LSd}
  into the mold of Figure~\ref{fig:controlsetup} as follows: consider a fictitious
  plant $G$ with input $u(t)$ and output $y(t)$ that is  determined by the state-space model
  \begin{eqnarray} \dot x(t) & = & - \nabla f\left( x(t) \right) + u(t) \label{eq:f1} \\ y(t) & = & L x(t). \label{eq:f2} \end{eqnarray}
In a sense, this fictitious plant models the agents selfishly minimizing their private cost function---an accurate statement when the input $u(t)$ is off duty. The role of $u(t)$ is thus to nudge the agents away from their greedy trajectories so as to build a consensus, that is, so as to have $L x(t) = 0$, which is the desired reference signal for the output $y(t)$. Thus, the error signal is $e(t) = -y(t)$. If we now adopt a proportional-integral type of controller $u(t) = e(t) + L \int_{-\infty}^t e(s) ds$ and plug it in the fictitious plant~\eqref{eq:f1}-\eqref{eq:f2}, we arrive at the primal-dual algorithm~\eqref{eq:LSp}-\eqref{eq:LSd}.

Interpreting an algorithm such as~\eqref{eq:LSp}-\eqref{eq:LSd}
as a certain kind of controller expands our understanding of the algorithm because many insights from control can be brought to bear, for example, insights about the qualitative behavior of the algorithm regarding classic control metrics such as overshoot or settling time.
A good case in point is the work~\cite{Hu2017},
 in which known properties of several standard first-order methods such as the gradient method,
 the Nesterov accelerated gradient method~\cite{Nesterov-Gradient}, and the heavy-ball method~\cite{PolyakHB},
 find pleasant explanations through the  angle of loop-shaping control synthesis.
   Furthermore, control-theoretic tools have been proved to be very
   useful in the analysis and design of primal-dual distributed methods.
   An excellent example is the small-gain theorem, e.g., \cite{SmallGainTheorem}, a standard tool for analyzing stability in the context of
interconnected dynamical systems, utilized, e.g., in~\cite{SmallGainNedicTimeVar,UnificationAndGeneral,SmallGainNedicUncoord}.
The tool has been proved to be well-suited, e.g.,
in analyzing the methods in time-varying networks settings,
 in which case it provides an elegant way to cope with the challenging time-varying dimension of the problem.
As another example, \cite{LessardDistributed} (see also~\cite{LessardIQC}) provides a method
  for analysis of distributed
   optimization algorithms, where a small-sized semidefinite program (SDP),
    with the size independent of the network and optimization variable dimensions,
    is associated to the algorithm of interest. Then,
    the numerical solution of the SDP gives a tight upper bound on
    the algorithm convergence rate. Furthermore,
    the authors provide a novel algorithm dubbed SVL,
    designed based on the utilized control-theoretic approach, i.e., through optimization of the associated SDP.
 The recent work~\cite{Sundararajan2019},
which advances a canonical form for distributed first-order methods that not only encapsulates many known state-of-art algorithms
but also explains, for example, why the early, slower single-state first-order variants could not do away with the diminishing stepsize.

As a final example, techniques from control also help in studying time-varying (convex) optimization problems, e.g.,~\cite{Simonetto}.
To illustrate, consider an optimization problem that depends on a parameter.
This parameter thus sets implicitly the solution of the problem,
and commands explicitly the iterations of any associated primal-dual algorithm.
Imagine now that the parameter starts changing stochastically, as fast as the algorithm iterations are taking place.
Such fast-changing solution cannot be tracked perfectly by the algorithm iterations, which raises the natural performance question:
how good is this tracking, say, as measured by the mean-squared-error (MSE)?
In control parlance, this question is about the L2 gain from an input  signal (the time-varying parameter)
to an output signal (the error between the iterations and the time-varying solution), with the underlying virtual
system being the algorithm. Such an issue was tackled in~\cite{ChenLau},
in the context of a primal-dual algorithm for wireless applications,
by resorting to Lyapunov and LaSalle's standard techniques.
As another example of this line of work,
the more recent reference~\cite{Porco} uses arguments from passivity theory to study the L2 gain of
a primal-dual algorithm from random disturbances to deviations around the problem solution.

\vspace{-4mm}
\section{Overview of recent results: Inexact updates, complexity, and tradeoffs}
\label{section-overview}
This Section surveys recent results on
ALM and primal-dual optimization methods for master-worker and
fully distributed architectures introduced in Section~{II}.
 A major focus recently has been on designing and analyzing methods with
 inexact primal updates. While this is a traditional topic,
 a more recently effort has been put on how inexact updates
 translate into complexity of the overall algorithm, measured, e.g.,
 in terms of the overall number of inner (primal) iterations to achieve
 a prescribed accuracy. This kind of analysis more closely
 reflects the algorithm performance in terms of computational and communication costs
 with respect to traditional metrics in terms of the number of dual
 (outer) iterations. Subsection~{IV-A}
  considers master-worker architectures,
  Subsection~{IV-B} studies fully distributed architectures,
  while Subsection~{IV-C} presents very recent results on
  lower complexity bounds for fully distributed architectures, i.e.,
  consensus optimization.

\subsection{Dual decomposition framework}
\label{subsection-dual-decomp}

The conventional (dual)
 decomposition framework for the augmented Lagrangian
 methods has received renewed interest in the past decade,
 e.g., \cite{NecoaraSykensSmoothing,NecoaraMPC2},
 and has been applied in various fields such
 as distributed control, e.g.,~\cite{JoaoMotaMPC,NecoaraMPC2}.

The recent results are mostly concerned with establishing complexity in terms of
the primal solution sub-optimality and
primal feasibility gaps. Therein,
the focus is to explicitly account for and model
the inexactness of the primal updates, which are for example carried out by a
specified (inner) iterative method.
Then, the goal is to establish complexity
 results in terms of the total number of inner and outer
 iterations for convergence.
 This translates into  precise
 quantification of the communication-computational tradeoffs
 for solving different classes of problems.

Reference~\cite{NecoaraMPC2}
considers a dual and fast dual gradient methods and assumes that
the primal problems are solved up to a
pre-specified accuracy $e_{\mathrm{in}}$.
The paper then establishes the complexity of the overall
algorithm as a function of $e_{\mathrm{in}}$
and discusses optimal setting of $e_{\mathrm{in}}$
 so that the overall complexity is minimized.
 Reference~\cite{NecoaraAUFLagr}
 considers a similar framework for \emph{augmented} Lagrangian methods
 and establishes complexity of the overall inexact augmented Lagrangian methods.
 For an inexact accelerated dual gradient method,
 the reference establishes the complexity
 (in terms of outer iterations) of order~$O\left( 1/\sqrt{\epsilon}\right)$,
 where the inner problems need to be solved up to accuracy~$O(\epsilon\sqrt{\epsilon})$.
  An interesting approach utilizes smoothing
  of the dual function~\cite{NesterovSmoothing}; we refer to~\cite{NecoaraSykensSmoothing}
  for the method and the related complexity results.
  Reference~\cite{GeorgeLanAL}
  considers the augmented Lagrangian method when
  the primal updates are carried out via an
  inexact Nesterov gradient method (inner iterations).
  The paper then establishes complexity
  of the method in terms of the overall number of inner iterations,
  and it also proposes a method with improved complexity
   based on the addition of a strongly convex term to the
   original objective function.


The revival of AHU-type methods (see \eqref{eq:ALprimalupdate-AHU-Sec-2}--\eqref{eq:ALdualupdate-AHU-Sec-2}) in the past one to two decades
(see, e.g., \cite{ChambolleImaging,EsserImaging,GeorgeLanPDold,GeorgeLanSaddlePoint,NemirovskiExtragradient,MonteiroSweiterExtragradient,NesterovSmoothing})
   is partly motivated by their success in imaging applications, e.g.,~\cite{ChambolleImaging,EsserImaging}.
  A canonical setting considered is on non-smooth convex-concave saddle point problems.
 Reference \cite{ChambolleImaging}
 considers a unified framework of AHU-type algorithms
 and shows that the methods converge at the optimal rate $O(1/k)$,
 where $k$ is the iteration counter.
 Another line of work~\cite{NemirovskiExtragradient,MonteiroSweiterExtragradient}
 considers primal-dual methods based
 on the extragradient method~\cite{ExtragradientKorpelevich}
 and demonstrates for methods therein also the rate of convergence~$O(1/k)$.
 Reference~\cite{NesterovSmoothing} shows that the smoothing
 technique due to Nesterov also allows for optimal
 rate~$O(1/k)$. While the works~\cite{ChambolleImaging,NemirovskiExtragradient,MonteiroSweiterExtragradient}
 do not exhibit optimal rates in terms of
  other problem parameters like the function's Lipschitz constants,
  the smoothing technique~\cite{NesterovSmoothing} coupled with
  the Nesterov gradient method exhibits optimal scaling with respect to
  problem parameters as well. Finally,
  reference~\cite{GeorgeLanSaddlePoint}
  constructs a novel method that also matches
  the optimal rate with respect to the iteration counter and problem parameters.

\subsection{Consensus optimization}
\label{subsection-overview-consensus}

We now provide an overview of
primal-dual algorithms for distributed consensus
optimization. A starting point for many of
the algorithms is~\eqref{eqn-AL-generic-ALM-primal}--\eqref{eqn-AL-generic-ALM-dual}.
 We classify here the algorithms with respect to which variant of the primal-dual method
 they use as a ``baseline'', i.e.,
 we divide them into the classical ALM-based, ADMM-based, or AHU-based methods.

 \textbf{ALM-based methods}.
 Several ALM-based methods have been proposed, e.g., \cite{cooperative-convex,LinearALRateDJJ,ErminWeiPDframework,NEWRibeiroPDQN,RibeiroESOM,FrensisBachDistrOptimal}.
 Reference~\cite{LinearRateAL} considers several deterministic and randomized variants
 to solve~\eqref{eq:ALp} inexactly in an iterative fashion,
 including gradient-like and Jacobi-like primal updates.
 Reference~\cite{ErminWeiPDframework} considers
 gradient-like primal updates and proposes
algorithms to reduce communication and computational costs with theoretical
guarantees. It also shows by simulation
 that performing a few inner iterations (2-4, more precisely)
  usually improves performance over the
  single inner iteration round-methods like EXTRA~\cite{Shi2015}.

Reference \cite{RibeiroESOM} considers an inexact
ALM based on
an augmented Lagrangian function with a proximal term added.
 The primal variable update is then replaced
 by a quadratic (second order) approximation
 of the cost function.
 Reference \cite{RibeiroESOM}
 shows that, provided that the dual step size
 is below a threshold value,
 the method exhibits a globally linear convergence rate.

  Very recently, \cite{NEWRibeiroPDQN}
  proposes a primal-dual method based on augmented Lagrangian
that employs a second order, quasi-Newton type update both on
  the primal and on the dual variables.

  \textbf{ADMM-based methods}.  An important class of primal-dual methods is based on ADMM, e.g.,
  \cite{ADMMlinear,RibeiroLinearized,DQM,LessardDistributed,Joao-Mota-3}.
Recent works approximate the objective function
 in the primal variable update
 via a linear or a quadratic approximation,
 \cite{RibeiroLinearized,DQM}.
  In this way, one can trade-off the communication and computational
  costs of the algorithm. References~\cite{RibeiroLinearized,DQM}
  show that, through such approximations, usually a significant
  computational cost reduction can be achieved
  at a minor-to-moderate additional communication cost.
   Recently, \cite{LingLinearizedCensored} proposes a censored
   version of the linearly approximated ADMM in \cite{RibeiroLinearized}
   to further improve communication efficiency
   and establishes its linear convergence rate guarantees.

  \textbf{AHU-based methods}. Several recent works have appeared which
 are based on AHU-type dynamics \eqref{eq:ALprimalupdate-AHU-Sec-2}--\eqref{eq:ALdualupdate-AHU-Sec-2},
 e.g., \cite{MartinezPrimalDual,NedicPerturbation,Shi2015,SmallGainNedicUncoord,UnificationAndGeneral,LessardDistributed,SayedExact1}.
  References~\cite{Martinez,YuanXuZhaoPrimalDual}
 are among early works on AHU type methods for solving \eqref{eq:dopt}.
 The methods consider very generic
 costs and constraints and propose methods that
 build upon the classical AHU-type dynamics~\eqref{eq:ALprimalupdate-AHU-Sec-2}--\eqref{eq:ALdualupdate-AHU-Sec-2}.
 Reference \cite{YuanXuZhaoPrimalDual} allows for
 public constraints (known by all agents), while reference \cite{Martinez}
 also allows for private constraints.
The work \cite{CortesSaddlePoint2014ContTime}
considers continuous time AHU dynamics and proposes a method
variant that converges for balanced directed graphs.
Reference \cite{CortesSaddlePoint2016}
considers more generic convex-concave saddle point problems that,
when applied to problem \eqref{eq:dopt},
translate into an AHU-type method.
Reference \cite{NedicPerturbation} proposes
a primal-dual method for more generic problems but which can be applied to~\eqref{eq:dopt},
based on the perturbed primal-dual AHU dynamics~\cite{PerturbedPrimalDual}.
 For the case of general convex-concave saddle-point problems, they establish convergence of the running time-averages of the local estimates to a saddle point under periodic connectivity of the communication digraphs.

References~\cite{SayedExact1,SayedExact2} develop
 exact variants of the popular diffusion methods, e.g., \cite{AliSayedIEEEproc,SayedEstimation}.
 A distinctive feature of these methods is that they relax the requirement
 that the underlying weight matrix be \emph{doubly} stochastic, i.e.,
 the matrix is allowed to be singly stochastic.
  The works~\cite{Usman1,Usman2,Usman3} devise
  exact first order distributed methods for directed networks.
 The work~\cite{AliSayedVazno}
 provides a very general framework of proximal exact methods which subsumes
 many existing primal-dual methods, including~\cite{SmallGainKin1,SmallGainNedicUncoord,harnessing,SmallGainNedicTimeVar,RibeiroLinearized,Shi2015,SayedExact1,NIDSweiShi}.

Table~I summarizes a number of primal-dual methods for consensus optimization, classified
according to which baseline method they use (classical ALM, ADMM, or AHU).
Table~I describes what type of primal and dual updates the algorithms utilize,
whether they are single or double looped, and whether they
are proposed for unconstrained or constrained optimization problems.

\begin{table*}[!ht]
\setlength\extrarowheight{2pt} 
\begin{tabularx}{\textwidth}{|X|X|X|X|X|}
\hline
\small{Baseline method}          &  \small{Primal update} & \small{Dual update} & \small{Presence of double loop}
& \small{Handling constraints}
\\ \hline
\small{ALM:  \cite{cooperative-convex,LinearALRateDJJ,ErminWeiPDframework,NEWRibeiroPDQN,RibeiroESOM,FrensisBachDistrOptimal}} &
\small{\cite{LinearALRateDJJ}: multiple gradient or Jacobi steps on AL;
 \cite{ErminWeiPDframework}: multiple gradient steps on AL; \cite{RibeiroESOM}:
  quadratic approximation of AL;
 \cite{NEWRibeiroPDQN}: quasi-Newton approx. of AL;
 \cite{FrensisBachDistrOptimal}: exact minimization of Lagrangian}
 &
 \small{\cite{LinearALRateDJJ,ErminWeiPDframework}: dual gradient ascent;
   \cite{NEWRibeiroPDQN}:
   quasi-Newton dual ascent;
 \cite{FrensisBachDistrOptimal}: accelerated (Nesterov) dual gradient ascent on the Lagrangian} & \small{
   \cite{cooperative-convex,LinearALRateDJJ,ErminWeiPDframework,FrensisBachDistrOptimal}: double loop;
    \cite{NEWRibeiroPDQN,RibeiroESOM}: single loop
   } & \small{\cite{cooperative-convex}: constrained problems;
   \cite{LinearALRateDJJ,ErminWeiPDframework,NEWRibeiroPDQN,RibeiroESOM,FrensisBachDistrOptimal}:
   unconstrained problems}
    \\
   \hline
\small{ADMM: \cite{ADMMlinear,RibeiroLinearized,DQM,LessardDistributed,Joao-Mota-3}} &
\small{\cite{ADMMlinear,Joao-Mota-3}: minimize AL exactly; \cite{RibeiroLinearized,LessardDistributed}: minimize linear approx. of AL;
\cite{DQM}: minimize quadratic approximation of AL} & \small{dual gradient ascent}  & \small{
single loop
} & \small{\cite{Joao-Mota-3}: constrained problems;
\cite{RibeiroLinearized,DQM,LessardDistributed,ADMMlinear}: unconstrained problems}
\\ \hline
\small{AHU:  \cite{MartinezPrimalDual,NedicPerturbation,Shi2015,SmallGainNedicUncoord,UnificationAndGeneral,LessardDistributed,SayedExact1}}
&  \small{single gradient step on AL} & \small{dual gradient ascent} & \small{single loop} & \small{
\cite{MartinezPrimalDual,NedicPerturbation}: constrained problems;
 \cite{Shi2015,SmallGainNedicUncoord,UnificationAndGeneral,SayedExact1,LessardDistributed}: unconstrained problems
} \\ \hline
\end{tabularx}
\caption{Overview of primal-dual methods for convex consensus optimization}
\end{table*}

\subsection{Optimal methods and lower complexity bounds for consensus optimization}
\label{section-methods-optimality-and-lower-complexity-bounds}
Very recently, lower complexity bounds for
consensus optimization have been analyzed, and
optimal algorithms that can achieve these lower bounds have been proposed.
 These results
 determine what is the lowest possible (appropriately measured) communication and computational cost
 for any distributed algorithm and also give insights into the interplay between the two costs.
 Interestingly, in certain cases the optimal algorithms and matching complexity bounds
 are derived based upon or admit ALM or primal-dual interpretations.

The work~\cite{FrensisBachDistrOptimal} introduces a
 formal black-box oracle model that allows
to precisely account for the complexity bounds.
 Therein, each node has a finite local memory, and it
 can evaluate the gradient and the Fenchel conjugates of its local cost function~$f_i$.
 Further, each node $i$ can, at time $k$, share a value with all or part of its neighbors.
 The class of admissible algorithms is then roughly restricted
 to the ones that generate the solution estimate as
 arbitrary linear combinations of the past gradients and Fenchel conjugates
 available in node $i$'s local memory. We refer to~\cite{FrensisBachDistrOptimal}
 for a more formal and detailed description of the oracle model.

  We denote by $\mathcal G$ the underlying communication graph. Next, quantity $\tau$ is the unit communication cost relative to a unit computational cost, i.e.,
 $\tau$ is the time to communicate a value to a neighbor in the network,
 while the time to perform a unit computation equals one.
 We consider an arbitrary algorithm that conforms to
 the black box model above
 and utilizes the averaging matrix $W$
 that has the spectral gap~$\gamma$.

Reference~\cite{FrensisBachDistrOptimal} shows that, for any $\gamma > 0$,
there exists a weight matrix~$W$ of spectral gap $\gamma$, and
 ${{{\mu}}}$-strongly convex, $\ell $-smooth functions
$f_i$'s such that, with $\kappa=\ell/{{{\mu}}}$, for any black-box procedure using~$W$ the time to reach a precision $\varepsilon > 0$ is lower bounded by
 $
\Omega\left(\sqrt{\kappa}\left(1 + \frac{\tau}{\sqrt{\gamma}}\right)\ln\left(\frac{1}{\varepsilon}\right)\right).
 $
Reference~\cite{FrensisBachDistrOptimal} further provides an optimal algorithm
that matches the lower complexity bound above.
 The algorithm requires a priori knowledge of network parameters
such as spectral gap~$\gamma$ and some knowledge on the eigenvalues of matrix~$W$.
 It is also interesting to compare
 the achieved convergence time with~\cite{NIDSweiShi}.
  Reference~\cite{NIDSweiShi}
  derives a convergence time for the algorithm therein
  of order $O\left( \max\{\kappa,\,1/\gamma \}\mathrm{ln}(1/\varepsilon)\right)$; the time exhibits a
  worse dependence on the condition number~$\kappa$ than in~\cite{FrensisBachDistrOptimal}, but on the other hand
    it does not depend on the product of $\gamma$ and
    $\kappa$, but on their maximum.

Reference~\cite{ScutariPrimalDualACC2019} considers a similar black box oracle model
but allows for lower ``oracle power,'' wherein the nodes' ability
to calculate Fenchel conjugates is replaced by the nodes'
ability to evaluate gradients of their local functions~$f_i$'s.
 The reference assumes
 that the functions $f_i$'s
 are convex and have Lipschitz continuous gradients with Lipschitz constant~$\ell$.
 The authors show that, for any $\gamma > 0$, there exists a weight matrix~$W$ of spectral gap $\gamma$ and
smooth functions $f_i$'s with $\ell$-Lipschitz continuous gradients such that,
for any black-box procedure based on the considered oracle using~$W$, the time to reach a precision $\varepsilon > 0$ is lower bounded by
 $
	\Omega{ \left( \left(1+\frac{\tau}{\sqrt{\gamma}}\right)   \sqrt{\frac{L}{\varepsilon} } \right)}.
 $

The work~\cite{FrensisBachDistributedOptNonsmooth}
 establishes lower complexity bounds for non-smooth functions $f_i$'s and a black box oracle model
similar to \cite{ScutariPrimalDualACC2019}.

We close the section by noting that, while some lower optimality bounds have been met,
 there is still room for improvements.
 Namely, the current methods that achieve the optimality bounds are often constructive and may be complicated,
 and may require an a priori knowledge of several global system parameters,
 which may not be available in practice.

\section{Applications}
\label{section-applications}
In this section, we present two applications of
ALM and primal-dual methods. In Subsection~{V-A},
we introduce for the first time ALM-type methods
in the context of federated learning; more precisely,
we describe how PDMM \eqref{eqn-rpdmm_xi}--\eqref{eqn-rpdmm_yhat}
 can be used for federated learning.
 Subsection~{V-B} considers distributed energy trading in microgrids,
 wherein we provide some novel insights.

\subsection{Federated learning}
\label{subsection-federated-learning}

Federated learning is an emerging
paradigm where a potentially large number of (edge) nodes, or
multiple clients (e.g., collaborative organizations),
cooperatively learn a global model while performing
only local update recommendations based on their local data. We refer
here to the edge nodes or clients as devices. With federated learning,
the devices
 keep their data locally, without the need for raw data communication,
 e.g.,~\cite{FederatedLearningIntro,FederatedLearningOverview,CrossSiloFederatedLearning}.

As pointed out in~\cite{FederatedLearningOverview},
the setting of federated learning poses several challenges to
 the design of learning and optimization algorithms; we highlight
 some of these challenges.
 First,
 due to the communication bottleneck, only a small subset
 of devices may be allowed to communicate with the server at a single round.
 Second, the devices may mutually significantly differ
 in their communication, computational, and storage capabilities.
 Third, the data available at different devices in general follows
 a different statistical distribution, where the distribution of
 the data is unknown both at the device and server levels.

More formally, the problem of interest is as follows. There are $N$ devices and a single server, where each device can communicate
with the server bidirectionally. The goal is to solve the following optimization problem:
 \begin{equation}
 \label{eqn-opt-prob-original-federated}
\mathrm{minimize}_{x \in {\mathbb R}^d}\,\,f(x) := \sum_{i=1}^{N} p_i\,F_i(x).
\end{equation}
Here,
$p_i > 0$, and $\sum_{i=1}^N p_i =1$.
We assume that $p_i$ is known by device
$i$ and that all $p_i$'s are known at the server.
 Furthermore, $F_i(x)$ is the empirical risk with respect to
the data available at device $i$. The distributions that generate the different devices' data can be different.
 Denote by
$f_i(x):=p_i\,F_i(x)$. We also let $$F(x) = F(x_1,....,x_N) = \sum_{i=1}^N f_i(x_i).$$

 Recently proposed algorithms to solve \eqref{eqn-opt-prob-original-federated}
  include, e.g., \fedavg~\cite{FederatedLearningIntro}, and \fedprox~\cite{FederatedLearningAnitHetero}.

We describe here how PDMM,~\cite{PDMMTomLuo},
 can be used to solve~\eqref{eqn-opt-prob-original-federated}.
 First, we introduce the following (equivalent) reformulation of \eqref{eqn-opt-prob-original-federated},
 analogous to the reformulation in~\cite{BoydADMM} for consensus optimization. That is, we
 minimize
 $
F(x)=\sum_{i=1}^N f_i(x_i)$, subject to the constraints
$z=x_i,$ $i=1,...N.$

 The augmented Lagrangian function is given by:
\[
\mathcal{L}(x,z;{{{\lambda}}}) = F(x) + \sum_{i=1}^N {{{\lambda}}}_i^\top \left( z - x_i\right)
+\frac{\rho}{2}\sum_{i=1}^N \|z-x_i\|^2.
\]

 We denote by $S_k \subset \{0,1,...,N\}$
 the variable blocks that
 the server selects to update via PDMM at iteration~$k$.
 Here, we assign index $i$ to variable block $x_i$, $i \in \{1,...,N\}$,
 and index $0$ to variable block~$z$.
 The blocks in $S_k$ can be selected uniformly at random from
 set $\{0,1,...,N\}$, such that the cardinality of $S_k$ equals~$M$.
We apply
 here the PDMM
variant without the backward dual step, i.e.,
we assume that the weights  that correspond
to the Bregman divergence term in~\eqref{eqn-rpdmm_xi} are sufficiently large.
 We let the Bregman divergence term in~\eqref{eqn-rpdmm_xi} associated with
 the $z$-variable be
 $\frac{\eta_0}{2}\|z-z^{(k)}\|^2$,
 while the Bregman divergence terms associated with
 the $x_i$'s are free to choose.
 Next, taking advantage of the fact that
 the $z$-update can be carried out in closed form,
 and simplifying further,
 one can arrive at the method
 presented in Algorithm~1.
 Variables $z,x_1,...,x_N$,
 can be, e.g., initialized to zero.
\begin{algorithm}[h]
	\begin{algorithmic}
		\caption{A PDMM-based method for federated learning}
		\label{algorithm-fedprox}	
		\STATE {\bf Input:}  $K$, $T$, $x^{(0)}$, $N$, $p_i$, $i=1,\cdots, N$
		\FOR  {$k=0, \cdots, T-1$}
			\STATE (S1) Server selects a subset $S_k$ of variable blocks
            \STATE (S2) Server sends $z^{(k)}$ to each device $i$, $i \in S_k$.
			\STATE (S3) Server sends the dual variable ${{{\lambda}}}_i^{(k)}$ to $i$, $i \in S_k$.
            \STATE (S4) If~$0\in S_k$, the server performs the following update:
			{\small$z^{(k+1)} = \frac{1}{\rho\,N+\eta_0}
\left( \rho \,\sum_{i =1}^N  x_i^{(k)} + \sum_{i=1}^N \lambda_i^{(k)} + \eta_0\,z^{(k)}\right)$;
			\STATE (S5) Each chosen device $i \in S_k$ finds $x_i^{(k+1)}$ as a  minimizer of:
			{\small$  f_i(x) + ({{{\lambda}}}_i^{(k)})^\top (x- z^{(k)}) +\frac{\rho}{2}\|x-z^{(k)}\|^2 +
\eta_{i}^{(k)} B_{i}(x, x_i^{(k)})$}\
			\STATE (S6) Each device $i \in S_k$ sends $x_i^{(k+1)}$ back to the server
            \STATE (S7) Each device $i \notin S_k$ sets locally $x_i^{(k+1)}=x_i^{(k)}$ (no communication required)
            otherwise, it sets $z^{(k+1)}=z^{(k)}$} 
            \STATE (S8) Server updates all dual variables ${{{\lambda}}}_i$, $i=1,...,N$, as follows:
            ${{{\lambda}}}_i^{(k+1)} = {{{\lambda}}}_i^{(k)} + \rho\, \left( x_i^{(k+1)} - z^{(k+1)} \right)$
		\ENDFOR
	\end{algorithmic}
\end{algorithm}
We now comment on Algorithm~1. Interestingly,
Algorithm~1
has some similarities with \fedprox in \cite{FederatedLearningAnitHetero},
even though it is derived from a very different perspective, but it also has the following
key differences. First, while in \fedprox device $i$ minimizes
(approximately) function $F_i(x)+\frac{{{{\rho}}}}{2}\| x-z^{(k)}\|^2$ (see~\cite{FederatedLearningAnitHetero}),
 with PDMM the latter function is replaced with\footnote{
 It can be shown~\cite{PDMMTomLuo} that the exact minimization of function
 $  f_i(x) + ({{{\lambda}}}_i^{(k)})^\top (x- z^{(k)}) +\frac{\rho}{2}\|x-z^{(k)}\|^2 +
\eta_{i}^{(k)} B_{i}(x, z^{(k)})$ in Algorithm~2
corresponds to an inexact minimization of function
$  f_i(x) + ({{{\lambda}}}_i^{(k)})^\top (x- z^{(k)}) +\frac{\rho}{2}\|x-z^{(k)}\|^2 $.
Hence we abstract the influence of the term $\eta_{i}^{(k)} B_{i}(x, z^{(k)})$.
}
 $$f_i(x)+({{{\lambda}}}_i^{(k)})^\top \left(  x - z^{(k)} \right) + \frac{\rho}{2}\| x-z^{(k)}\|^2.$$
   The difference in the inclusion of $f_i$ with PDMM
  instead of $F_i$ is minor\footnote{Actually,
  ignoring the term $({{{\lambda}}}_i^{(k)})^\top \left( z^{(k)} - x \right) $, and
  allowing for different penalty parameters $\rho_i$ across
  different nodes $i$, one can divide
  the function $f_i(x)+ \frac{\rho}{2}\| x-z^{(k)}\|^2$
  by $p_i$, $p_i>0$, to arrive at the function used by \fedprox.}.
  However, the difference in accounting for
  the dual variable ${{{\lambda}}}_i^{(k)}$ is important.
  Intuitively, both the proximal term $\frac{{{{\lambda}}}}{2}\| x-z^{(k)}\|^2$
  with \fedprox and the term ${{{\lambda}}}_i^{(k)})^\top \left(  x  - z^{(k)} \right) + \frac{\rho}{2}\| x-z^{(k)}\|^2$
   with PDMM serves to combat  the drift that may appear due to the
   difference between the local function $F_i(\cdot)$ and the global function
   $\sum_{i=1}^N p_i\,F_i(\cdot)$.
   However, the effect may be strengthened with PDMM, due to the inclusion of the dual
   variable. Second,
   the two algorithms utilize different aggregation updates -- see
   (S4) in Algorithm~1.
   That is, the aggregation step in \fedprox  takes into
   account only the $x_i$'s for $i \in S_k$, i.e., the $x_i$'s of active devices only.
    On the other hand,
    PDMM accounts
   for the latest historical value of $x_i$ from each device~$i$,
   even when at iteration~$k$ device $i$ is not active.
   This may ``smoothen'' the update of $z$ and implicitly
   incorporates history in the $z$-update rule.

When
each of the functions~$F_i$ is convex,
 PDMM enjoys global convergence,
 and the expected optimality gap
  $\mathbb{E}[f(z^{(k)})-f^\star]$ converges to zero,
  where $f^\star$ is the optimal value of~\eqref{eqn-opt-prob-original-federated}.

In the context of federated learning, it is also highly relevant to consider non-convex problems,
 e.g., for training neural networks.
 Interestingly, we can apply
 the results of~\cite{NonconvexADMM}
  to obtain an algorithm similar to
   Algorithm~1 that works for non-convex $F_i$'s as well.
   More precisely, a direct application of
   Algorithm~2 in~\cite{NonconvexADMM} leads to
   an algorithm very similar to Algorithm~1;
   the only differences with respect to Algorithm~1 are the following:
   set $\eta_i^{(k)}=0$, for all $i,k$,
   $\eta_0=0$; in (S5) use $z^{(k+1)}$ instead of $z^{(k)}$; and (S8) is applied only
   to the $\lambda_i$'s for $i \in S_k$,
   otherwise they are kept unchanged in the current iteration.
    Applying the results in~\cite{NonconvexADMM},
    the resulting algorithm converges to a stationary solution
     almost surely, provided that  the $F_i$'s
      have Lipschitz continuous gradient, the penalty
      parameter~$\rho$ is large enough, and $f$ is bounded from below.
      (See Theorem 2.4 in~\cite{NonconvexADMM}.)


\subsection{Microgrid peer-to-peer energy trading}
We now present how augmented Lagrangian methods
can be applied to distributed
energy trading, e.g., \cite{MicroGridsKelvinJournal}.

Distributed energy trading has received
significant attention recently, e.g.,
\cite{MicroGridsKelvinJournal,MicroGridsKelvin1,MicroGridsKelvin2,MicroGridsLowVoltage,MicroGridsMatamorosJournal,MicroGridsNEW,MicroGridsPoor}.
 A major motivation for such studies comes from renewable energy generators.
 Many of these generators, called prosumers,
 can produce, consume, and trade energy with their neighboring (peer) prosumers
 within a trading graph, e.g., \cite{MicroGridsPoor}.

 More formally, we consider $N$ energy trading peers connected within an undirected graph $\mathcal{G} = (V, {\mathcal E})$,
 where $V$ is the set of $N$ peers and ${\mathcal E}$ is the set of edges designating
 the trading capabilities of the peers.\footnote{
More generally, $\mathcal{G}$ can be directed, but for simplicity of presentation we keep $\mathcal G$ undirected.}
In other words,
  the presence of arc $(i,j)$ means that
  peer $j$ can buy energy from peer $i$. (We assume that the trading is symmetric, i.e.,
  if $(i,j) \in {\mathcal E}$, then $(j,i) \in {\mathcal E}$ as well.)
  We denote by $\Omega_i = \{j : (j,i) \in {\mathcal E}\}$ the neighborhood of $i$, i.e.,
the set of peers that are allowed to sell energy to peer $i$.

At
each time interval, peer $i$ generates $E_i^{(g)}$
and consumes $E_i^{(c)}$ energy units.
 If $(i,j) \in {\mathcal E}$, then peer $i$ can sell energy to peer $j$,
 in the amount of $E_{ij}$. We assume that the peers do not
 have energy storage capacity,
 and therefore each peer obeys the
 energy balance equation, as follows:
 \begin{equation}
 \label{eqn-energy-conservation}
 E_i^{(g)} + \sum_{j \in \Omega_i} E_{ji}
 = E_i^{(c)} + \sum_{j \in \Omega_i} E_{ij}.
 \end{equation}
Denote by $C_i(E)$ the cost incurred at peer $i$  for
generating $E$ energy units, and by
$\gamma_{ij} (E)$ the cost of trading and transferring the amount of energy
$E$ from peer $i$
to peer $j$. Then, we can formulate
an optimization problem of finding
optimal energy trading transactions $E_{ij}$, $(i,j) \in {\mathcal E}$, such  that
the total system costs of generation and trading are minimized:

{\allowdisplaybreaks{
{\small{
\begin{eqnarray}
&\,&\underset{{E_{ij}:\,(i,j)\in {\mathcal E}}}{\mathrm{minimize}}\,\,\sum_{i=1}^N C_i \left( E_i^{(c)} + \sum_{j \in \Omega_i} (E_{ij}-E_{ji}) \right) \nonumber\\
\label{eqn-formulation-energy-trading}
&\,&
\phantom{\mathrm{minimize}_{E_{ij}:\,(i,j)\in E}\,\,}+\sum_{j \in \Omega_i} \gamma_{ji}(E_{ji}) \\
&\,&\mathrm{subject\,to}\,\,\,\,\,\,E_{ij} \geq 0, \,\,(i,j)\in {\mathcal E} \nonumber \\
&\,& \phantom{\mathrm{subject\,to}}\,\,\,\,\,\,E_i^{(c)} + \sum_{j \in \Omega_i} (E_{ij}-E_{ji}) \geq 0, \,\,i=1,...,N. \nonumber
\end{eqnarray}}}}}

Reference~\cite{MicroGridsMatamorosJournal}
proposes a dual decomposition algorithm to solve \eqref{eqn-formulation-energy-trading}, see Algorithm~2.
 The algorithm is derived based on a reformulation of \eqref{eqn-formulation-energy-trading}
 that introduces an additional variable $\epsilon_i$ and constraint per peer, $\epsilon_i =\sum_{j \in \Omega_i} E_{ij}$,
 and then dualizes these constraints, so that the resulting
 Lagrangian $\mathcal{L}$ becomes:
   \begin{eqnarray}
  &\,&\mathcal{L}\left( E_{ij}, \,(i,j) \in E;\,\lambda_i, \,i=1,...,N\right)
  = \nonumber \\
  &\,&\sum_{i=1}^N C_i ( E_i^{(c)}
  + \epsilon_i -\sum_{j \in \Omega_i}-E_{ji}) +
  \sum_{j \in \Omega_i} \gamma_{ji}(E_{ji})   \label{eqn-lagrangian-microgrids}\\
  &+& \sum_{i=1}^N \lambda_i\,\left( \epsilon_i - \sum_{j \in \Omega_i} E_{ij}\right) . \nonumber
  \end{eqnarray}
  In Algorithm~2,
$\mathcal{E}_i = \{\epsilon_i,\,E_{ji},\,j \in \Omega_i:\,
\epsilon_i \geq 0, \,E_{ji} \geq 0,\,j \in \Omega_i,\,E_i^{(c)} + \sum_{j \in \Omega_i} (E_{ij}-E_{ji}) \geq 0\}$,
and
$\mathcal{L}_i$ is given by:
  \begin{eqnarray}
    \label{eqn-lagrangian-microgrids-V3}
 \lefteqn{\mathcal{L}_i \left(\epsilon_i,\,E_{ji},\,j\in \Omega_i;\,\lambda_j,\,j \in \Omega_i \cup \{i\}\right) = } \nonumber \\
& &   C_i \left( E_i^{(c)} + \epsilon_i -\sum_{j \in \Omega_i}E_{ji}\right) +
  \sum_{j \in \Omega_i}\gamma_{ji}(E_{ji}) \nonumber\\
  &+& \lambda_i\,\epsilon_i - \sum_{j \in \Omega_i} \lambda_j \,E_{ji}.
   \nonumber
  \end{eqnarray}
  Note that Algorithm~2 is fully distributed, i.e.,
only peer-to-peer communication is utilized over iterations. In (S4) of Algorithm~3, $\alpha_k$ is the dual step-size. See also~\cite{MicroGridsMatamorosJournal} for an
interesting economics interpretation of
Algorithm~2.
\begin{algorithm}
\caption{Dual decomposition~\cite{MicroGridsMatamorosJournal}}\label{alg:distributed}
\begin{algorithmic}
\STATE Peer $i$ initializes $\lambda_i^{(0)}$.
    \STATE (S1) Peer $i$ sends $\lambda_i^{(k)}$,  to its neighbors in $\Omega_i$,
    and it receives $\lambda_j^{(k)}$, $j \in \Omega_i$.
    \STATE (S2) Peer $i$ computes $\epsilon_i^{(k)}$ and $E_{ji}^{(k)}$,
    $j \in \Omega_i$ by
      finding a minimizer of the following problem:
\begin{eqnarray}
\label{eqn-formulated-microgrids-subproblem}
&\,&\mathrm{minimize}\,\,\mathcal{L}_i \left(\epsilon_i,\,E_{ji},\,j\in \Omega_i;\,\lambda^{(k)}\right) \\
&\,&\mathrm{subject\,to}\,\,\,\,\, 
(\epsilon_i,\,E_{ji},\,j \in \Omega_i) \in \mathcal{E}_i. \nonumber
\end{eqnarray}
    \STATE (S3) Peer-$i$ informs peers-$j$, $j \in \Omega_i$,  about the energy it is willing
      to buy, namely $E_{ji}^{(k)}$, at the given price $\lambda_{j}^{(k)}$
    \STATE (S4) with energy requests $E_{ij}^{(k)}$ from neighboring peers, peer $i$
      computes
      \begin{equation}\label{eqn-clearing-market}
        \lambda_i^{(k+1)} \gets \lambda_i^{(k)} + \alpha_k\,
      \left(\epsilon_i^{(k)} - \sum_{j \in \Omega_i}E_{ij}^{(k)}\right).
      \end{equation}
    \STATE $k\gets k+1$ and go to (S1).
\end{algorithmic}
\end{algorithm}

While Algorithm~2 has an appealing structure,
convergence guarantees it provides may be
 weak for certain models of energy generation costs $C_i(\cdot)$'s and
 energy transfer costs $\gamma_{ij}(\cdot)$'s, for instance
 when not all of them are strictly convex functions.
 Such scenarios are highly relevant; in fact,
 linear functions are frequently used, e.g.,
  \cite{MicroGridsKelvin1,MicroGridsKelvin2}.
  In the context of distributed energy trading, we provide here
  novel insights how Algorithm~2 can be
  made amenable for more generic structures of $C_i(\cdot)$'s and $\gamma_{ij}(\cdot)$'s,
  applying the recent results on (inexact) ALM~\cite{NecoaraAUFLagr}.
  Namely, we construct the \emph{augmented} Lagrangian function as:
  {\small{
  \begin{equation*}
  \mathcal{L}_{\rho}\left( \{\epsilon_i\}, \,\{E_{ij}\};\,\lambda\right)
   = \mathcal{L}\left( \{\epsilon_i\}, \,\{E_{ij}\};\,\lambda\right) +
   \frac{\rho}{2}\sum_{i=1}^N \|\epsilon_i - \sum_{j \in \Omega_i} E_{ij}\|^2,
  \end{equation*}}}
   where $\mathcal{L}$ is given in \eqref{eqn-lagrangian-microgrids}.
   We then propose to
   utilize an inexact ALM as in \cite{NecoaraAUFLagr}.
   This modifies Algorithm~2 as follows.
    Step~(S2) in Algorithm~2 (which was done by the peers in parallel) is
    replaced with the following step:
\begin{eqnarray}
\label{eqn-formulated-microgrids-subproblem-V2}
&\,&\mathrm{minimize}\,\,\mathcal{L}_{\rho} \left(\{\epsilon\},\,\{E_{ji}\};\,\lambda^{(k)}\right) \\
&\,&\mathrm{subject\,to}\,\,\,\,\, 
(\epsilon_i,\,E_{ji},\,j \in \Omega_i) \in \mathcal{E}_i,\,i=1,...,N. \nonumber
\end{eqnarray}
Problem~\eqref{eqn-formulated-microgrids-subproblem-V2}
cannot be solved ``in one shot'' in a distributed fashion.
However, it can be solved iteratively, e.g., through
a standard gradient descent, over inner iterations
$r=1,...,R$. It is easy to see that the inner iterations
involve only peer-to-peer communications. Namely,
every peer $i$ sends current value of $E_{ji}$ to peer $j$
 and receives  the current value of $E_{ij}$ from peer $j$, $j \in \Omega_i$.
  Utilizing the results of~\cite{NecoaraAUFLagr}, we have that
    the overall (outer) problem
   is solved
    with accuracy $\varepsilon_{\mathrm{out}}$
    in~$O(1/\varepsilon_{\mathrm{out}})$
     outer iterations~$k$,
     where the inner problems (over iterations~$r$)
      need to be solved up to accuracy~$O(\varepsilon_{\mathrm{out}})$.
       Omitting details, we
       also point out that an improved complexity
       can be achieved when, instead of
       ordinary gradient descent, the Nesterov gradient descent is utilized.
       Also, it is possible to a priori determine the number
       of inner iterations~$R$ that guarantee that, at every outer iteration~$k$,
       the inner problem is solved up to the desired accuracy~$O(\varepsilon_{\mathrm{out}})$.
       It is worth noting that ALM has been considered before
       in the context of a related power systems problem~\cite{ALclassicalAPPpowerSystems}. However, this work is not concerned
       with the (inner) iteration complexity nor with a distributed, peer-to-peer
       structure, of the underlying computational system.


\vspace{-4mm}
\section{Conclusion}
\label{section-conclusion}
The paper provided a study of recent results on
augmented Lagrangian methods (ALM) and related primal-dual methods
for large-scale and distributed optimization.
We gave a gentle introduction to ALM and its variants
 and introduced control-theoretic tools for its analysis.
 The paper gave an overview of recent results on the topic,
 focusing on inexact updates, iteration complexity,
 and distributed consensus optimization.
 We capitalized on recent results in the field to
 draw novel insights in the context of two emerging applications:
 federated learning and smart grid.

\vspace{-4mm}

\bibliographystyle{IEEEtran}
\bibliography{IEEEabrv,bibliographyIEEEProceedings}

\end{document}